\documentclass[12pt,a4paper]{article}
\usepackage{times}
\usepackage{pdfsync,hyperref}

\usepackage{amssymb}
\usepackage{latexsym}
\usepackage{amsmath}
\usepackage{amscd}
\usepackage{theorem}
\usepackage{xy}
\xyoption{all}
\CompileMatrices

\usepackage{graphicx}
\usepackage{calc}

\newlength{\wcwidth}
\newlength{\wcheight}
\newcommand{\widecheck}[1]{\ensuremath{
\settowidth{\wcwidth}{#1}
\settoheight{\wcheight}{#1}
\addtolength{\wcheight}{1pt}
\makebox[0cm][l]{%
\raisebox{\depth+\wcheight}[0cm][0cm]{%
\scalebox{-1}{$\widehat{\hphantom{#1}}$}}}#1
\rule{0pt}{\wcheight+2.5pt}}}

\def\dquote[#1]{{\it #1}}

\normalsize

\makeatletter

\def\@@seccntfont{\bfseries\slshape}
\def\@@secheaderfont{\bfseries\upshape}
\def\@@precnt{{\upshape}}
\def\@@postcnt{{\upshape}}
\def\@@startsection#1#2#3#4#5#6{\if@noskipsec \leavevmode \fi\par
  \@tempskipa #4\relax\@afterindenttrue
  \ifdim \@tempskipa <\z@\@tempskipa -\@tempskipa
  \@afterindentfalse\fi
  \if@nobreak\everypar{}\else\addpenalty\@secpenalty\addvspace\@tempskipa\fi
  \@ifstar{\@dblarg{\@sect{#1}{#2}{#3}{#4}{#5}{#6}}}%
    {\@dblarg{\@@sect{#1}{#2}{#3}{#4}{#5}{#6}}}}
\def\@@sect#1#2#3#4#5#6[#7]#8{\ifnum
  #2>\c@secnumdepth\let\@svsec\@empty
  \else\refstepcounter{#1}\protected@edef\@svsec{\@seccntformat{#1}\relax}\fi\@tempskipa #5\relax
  \ifdim
  \@tempskipa>\z@\begingroup#6{\@hangfrom{\hskip#3\relax\@svsec}\interlinepenalty 
  \@M #8\@@par}\endgroup\csname #1mark\endcsname{#7}\else
  \def\@svsechd{#6{\hskip #3\relax\@svsec #8}\csname
    #1mark\endcsname{#7}}\fi\@xsect{#5}}
\def\setseccntfmt{\renewcommand{\@seccntformat}[1]{\S
  \csname the##1\endcsname\hspace{1ex}}}
\def\@@seccntfmt{\renewcommand{\@seccntformat}[1]{%
  {\@@seccntfont\@@precnt\csname
  the##1\endcsname}\@@postcnt\hspace{1ex}}}
\newcommand{\@@secnostar}[1][]{\def\tmpa{}\def\tmpb{#1}%
  \ifx\tmpa\tmpb\def\tmpa{-1ex}\else\def\tmpa{-\labelsep}\fi\@@seccntfmt\@@startsection
  {\@@name}{\@@level}{0mm}{-\baselineskip}{\tmpa}{\@@secheaderfont}{#1}}
\newcommand{\@@secstar}[2][]{\def\tmpa{}\def\tmpb{#1}%
  \ifx\tmpa\tmpb\def\tmpa{#2}\else\def\tmpa{#1}\fi\@@seccntfmt\@startsection
  {\@@name}{\@@level}{0mm}{-\baselineskip}{-\labelsep}{\@@secheaderfont}[\tmpa]{#2}}
\def\@@section{\def\@@name{section}\def\@@level{1}\@ifstar{\@@secstar}{\@@secnostar}}
\def\@@subsection{\def\@@name{subsection}\def\@@level{2}\@ifstar{\@@secstar}{\@@secnostar}}
\def\@@subsubsection{\def\@@name{subsubsection}\def\@@level{3}\@ifstar{\@@secstar}{\@@secnostar}}
\def\@@paragraph{\def\@@name{paragraph}\def\@@level{4}\@ifstar{\@@secstar}{\@@secnostar}}
\def\@@subparagraph{\def\@@name{subparagraph}\def\@@level{5}\@ifstar{\@@secstar}{\@@secnostar}}
\let\@@latexsection=\section\let\mysection=\@@section
\let\@@latexsubsection=\subsection\let\mysubsection=\@@subsection
\let\@@latexsubsubsection=\subsubsection\let\mysubsubsection=\@@subsubsection
\let\@@latexparagraph=\paragraph\let\myparagraph=\@@paragraph
\let\@@latexsubparagraph=\subparagraph\let\mysubparagraph=\@@subparagraph
\renewcommand\section{\@startsection{section}{1}{\z@}%
{-3.5ex\@plus-1ex\@minus -.2ex}%
{2.3ex \@plus .2ex}%
{\normalfont\bfseries}}
\renewcommand\subsection{\@startsection{subsection}{2}{\z@}%
{-2.0ex\@plus-1ex\@minus -.2ex}%
{1.5ex \@plus .2ex}%
{\normalfont\it}}
\newif\if@@lastcharstar
\def\@@xxlastcharstar#1{\gdef\@@prevchar{}\@@lastcharstarfalse\@@yylastcharstar#1\end}
\def\@@yylastcharstar#1{\ifx#1\end\def\tmpa{*}\ifx\@@prevchar\tmpa\@@lastcharstartrue\fi
  \let\@@next=\relax\else\def\@@prevchar{#1}\let\@@next=\@@yylastcharstar\fi\@@next}
\gdef\@thm#1#2{\@@xxlastcharstar{#1}\if@@lastcharstar\else\refstepcounter{#1}\fi
  \trivlist\@topsep \theorempreskipamount\@topsepadd
  \theorempostskipamount
  \@ifnextchar [{\@ythm{#1}{#2}}{\@begintheorem{#2}{\csname
  the#1\endcsname}\ignorespaces}}
\gdef\th@nonumplain{\normalfont\itshape
  \def\@begintheorem##1##2{\item[\hskip\labelsep\theorem@headerfont
  ##1]}%
  \def\@opargbegintheorem##1##2##3{%
    \item[\hskip\labelsep \theorem@headerfont ##1\ ##3]}}
\gdef\th@change{
  \def\@begintheorem##1##2{\item[\hskip\labelsep{\@@seccntfont
  \@@precnt##2\@@postcnt\hskip 1ex}\theorem@headerfont ##1]}%
  \def\@opargbegintheorem##1##2##3{\item[\hskip\labelsep{\@@seccntfont
  \@@precnt##2\@@postcnt\hskip 1ex}\theorem@headerfont ##1\ ##3]}}
\makeatother   
\def\newtheoremset[#1][#2]{{\theoremstyle{change}\newtheorem{-#1}[section]{#2}%
  \newtheorem{#1}[subsection]{#2}\newtheorem{*#1}[subsubsection]{#2}%
  \newtheorem{**#1}[paragraph]{#2}\newtheorem{***#1}[subparagraph]{#2}}
  {\theoremstyle{nonumplain}\newtheorem{#1*}{#2}}}
\theoremheaderfont{\normalfont\bfseries}
{\theorembodyfont{\normalfont\itshape}
\newtheoremset[theorem][Theorem]\newtheoremset[proposition][Proposition]
\newtheoremset[lemma][Lemma]\newtheoremset[corollary][Corollary]
\newtheoremset[conjecture][Conjecture]}\newtheoremset[sublemma][Sublemma]
{\theorembodyfont{\normalfont\rmfamily}
\newtheoremset[definition][Definition]\newtheoremset[remark][Remark]
\newtheoremset[notation][Notation]\newtheoremset[question][Question]}

\newcommand{\op}{\operatorname}

\newcommand{\Hom}{\op{Hom}}

\usepackage{hyperref}

\def\poorsc#1#2 {#1{\boldx\uppercase{#2}} }

\begin{document}
\title{Floer theoretic invariants for 3- and 4-manifolds}

\author{Yi-Jen Lee\\ {\small Institute of Mathematical Sciences,}
  \\{\small the Chinese
 University of Hong Kong}\\{\small Shatin, N.T. Hong Kong}}
 
\date{\small October 2014, revised June 2016}
\maketitle

\begin{abstract}Seiberg-Witten (Floer) theory, Ozsvath-Szabo's
  Heegaard Floer theory, Hutchings's embedded contact homology, in
  different stages of development, define (or are expected to define)
  packages of invariants for 3- and 4-manifolds (including manifolds
  with boundary and manifolds with certain types of corners). We
  describe what are known about their relationship, what are expected,
  and raise some questions along the way.\end{abstract}

\section{Floer homologies for 3-manifolds: \(HF\), \(HM\), and \(ECH\)}

Let \((M, \mathfrak{s})\) be a closed spin-c manifold. Three
Floer-theoretic invariants
(of differential structure) are associated to \((M, \mathfrak{s})\): the Heegaard Floer homology of
Ozsvath-Szabo,  
\(HF^+(M,  \mathfrak{s})\), the Seiberg-Witten-Floer homology (aka the monopole Floer
homology) of Kronheimer-Mrowka, \(\widecheck{HM}(M,  \mathfrak{s}, c_b) \),
and, given in addition a contact structure \(\xi\)  on \(M\), the
embedded contact homology of Hutchings-Taubes, \(ECH(-M,  [\xi,
\mathfrak{s}])\). Here, \([\xi,\mathfrak{s}]\) denotes an element of
\(H_1(M; \mathbb{Z})\) determined by the pair \((\xi,\mathfrak{s})\)
described in the following manner: By definition, a contact structure
on \(M\) is an oriented 2-plane field. Meanwhile, on an oriented
3-manifold, a spin-c structure can be identified with an equivalence
class of oriented 2-plane fields. (See e.g. the last section of
\cite{HL1} or \cite{HL2}). In view of this, a contact structure on
\(M\) determines a spin-c structure \(\mathfrak{s}_\xi\). Recall that
the set of spin-c structures on \(M\) is a torsor over \(H^2(M;
\mathbb{Z})\simeq H_1(M; \mathbb{Z})\). By assigning \(0\in
H_1(M;\mathbb{Z})\) to the spin-c structure \(\mathfrak{s}_\xi\), one
defines an isomorphism from \(H_1(M; \mathbb{Z})\) to the set of
spin-c structures that intertwines with the \(H^2(M; \mathbb{Z})\)
action on both sides.

 It is now known that these three Floer homologies are all
isomorphic: the first equivalence, \(\widecheck{HM}(M,  \mathfrak{s})\simeq ECH(-M,  [\xi,
\mathfrak{s}])\), is established by Taubes \cite{T:ECH}; the second
equivalence, \(\widecheck{HM}(M,  \mathfrak{s})\simeq HF^+(M, 
\mathfrak{s})\), by Kutluhan-Lee-Taubes \cite{KLT}; the third
equivalence, \(HF^+(M, 
\frak{s})\simeq ECH(-M, [\xi,\mathfrak{s}])\), by Colin-Ghiggini-Honda
\cite{CGH}. These invariants come equipped with rich algebraic
structures: they take value in a cyclically-graded module over the
graded ring \(\mathbb{A}_\dag(M):=\mathbb{Z}[U]\otimes
\bigwedge^*H_1(M;\mathbb{Z})/\text{Tor}\), where \(U\) is of degree
\(-2\) and elements in \(H_1(M;\mathbb{Z})/\text{Tor}\) are of degree
\(-1\).

Among the three Floer theories, the first two carry the following
additional features:  The Heegaard Floer \(HF^+\) and the Seiberg-Witten
\(\widecheck{HM}\) are respectively one of a system of four flavors
of Floer homologies in either theory. These four flavors are also modules over
\(\mathbb{A}_\dag(M)\), and they are related by two fundamental long exact
sequences. In the Heegaard Floer case, the four flavors are denoted
\(HF^-\), \(HF^\infty \), \(HF^+\) and \(\widehat{HF}\), and the
fundamental sequences are: 
\begin{equation}\label{eq:fundsq}\begin{split}
\cdots \to HF^-(M, \mathfrak{s})\to HF^\infty(M, \mathfrak{s})\to HF^+(M,
\mathfrak{s})\to \cdots\\
\cdots\to HF^-(M, \mathfrak{s})\to HF^-(M, \mathfrak{s})\to \widehat{HF}(M,
\mathfrak{s})\to \cdots
\end{split}
\end{equation}
The map \(HF^-(M, \mathfrak{s})\to HF^-(M, \mathfrak{s})\) in the
second sequence above is the \(U\)-action; thus the second exact
sequence is essentially the defining sequence of the flavor \(\widehat{HF}(M,
\mathfrak{s})\) from the \(\mathbb{Z}[U]\)-module structure of
\(HF^-\). The pair \(HF^-\) and \(HF^+\) satisfies certain duality
property.

In parallel, the four flavors of Seiberg-Witten-Floer homologies
corresponding to \(HF^-\), \(HF^\infty\), \(HF^+\), \(\widehat{HF}\)
are respectively denoted \(\widehat{HM}\), \(\overline{HM}\),
\(\widecheck{HM}\), and \(\widetilde{HM}\), and they satisfy the same
fundamental long exact sequences and duality property as the
Heegaard Floer homology described above. 

The aforementioned equivalence theorems are useful mostly due of the very
different geometric origins of the three Floer homologies. These theorems will
be stated more precisely in next section. In preparation, here we give a
minimal sketch of the three Floer homologies' respective 
setup and background. For some representative applications of these equivalence
theorems, see Section 4 below. 

The Seiberg-Witten invariants for closed 4-manifolds were discovered
and underwent rapidly development during the second half of the
1990's. Their relative simplicity in comparison to its predecessor, the Yang-Mills
theory, led to significantly shortened proofs of many major theorems
originally obtained via the
latter, as well as important new results. See e.g. \cite{D2} and
\cite{Ts} for surveys. Most notably, Taubes was able to establish a 
surprising equivalence of the Seiberg-Witten invariant and a version of Gromov
invariant for closed 4-manifolds \cite{T}, which have very 
different constructions. (This is referred to as ``Taubes's \(SW=Gr\)
theorem'' below). Very roughly speaking, the Seiberg-Witten invariant
of a closed 4-manifold \(X\) is
defined by counting (equivalence classes) of solutions to the Seiberg-Witten
equation on \(X\), while the Gromov invariant is defined by counting pseudo-holomorphic curves
in \(X\) when \(X\) is equipped with a symplectic form
\(\omega\). Recall that given a spin-c structure \(\mathfrak{s}\) on \(X\), the
Seiberg-Witten equation takes the following form:
\begin{equation}\label{eq:SW}
F_{A}^{+} -(\Psi^{\dag}\tau\Psi- i \varpi)=0
\quad \text{and}\quad \mathcal{D}_{A}^+\Psi=0,
\end{equation}
where \(A\) is a Hermitian
connection on the line bundle \(\det\, (\mathbb{S}^{+})\),
\(\mathbb{S}=\mathbb{S}^+\oplus S^-\) being the spinor bundle associated to
\(\mathfrak{s}\), and \(\Psi\) is a section of \(\mathbb{S}^{+}\).  The term
\(\Psi^{\dag}\tau\Psi\) stands for the bilinear
map from $\mathbb{S}^{+}$ to $i\Lambda^{+}$ that is defined using the Clifford
multiplication, and $\mathcal{D}_{A}^+: \Gamma
(\mathbb{S}^+)\to \Gamma(\mathbb{S}^-)$ and $\mathcal{D}_A^-: \Gamma
(\mathbb{S}^-)\to \Gamma(\mathbb{S}^+)$ are the 4-dimensional Dirac operators on \(X\) defined
by the metric and the chosen connection \(A\).  Lastly, 
\(\varpi\) is a self-dual 2-form on \(X\). It is often referred to as
the ``perturbation form'' in the equation. When \((X, \omega)\) is
symplectic, one may choose the metric \(g\) on \(X\) to so that \(\omega\)
is self-dual (and hence harmonic). This defines an almost complex
structure \(J\) so that \(g(\cdot, \cdot )=\omega (\cdot, J\cdot)\). Taubes considered perturbation
forms of the type \begin{equation}\label{eq:-large_r_pert}
\varpi=2r\omega-iF_{A_K}^+,
\end{equation}
where \(r\in
\mathbb{R}^+\), and \(A_K\) is the induced connection on the canonical line
bundle \(K\) associated to \(J\). Use Clifford multiplication by
\(\omega\) to split \(\mathbb{S}^+\) as a sum of eigen-bundles
\(E\oplus E\otimes K^{-1}\), and let \(\alpha\) denote the
\(E\)-component of \(\Psi\) under this decomposition. It was shown that as one takes
\(r\to \infty\), the zero locus of \(\alpha\) approaches,  in certain
technical sense,  a (possibly
disconnected) \(J\)-holomorphic curve in a homology class determined
by \(\mathfrak{s}\). 

If instead of closed 4-manifolds, one considers gauge-theoretic
equations on 4-dimensional cylinders \(\mathbb{R}\times M\),  one may
in principle construct an associated Floer-homology for 3-manifolds. 
The Seiberg-Witten version of this, \(HM\), while long expected,
has rather technical actual construction and was not fully
written down until almost a decade later \cite{KM}. Nevertheless, in
comparison to the Yang-Mills version, usually called the
instantion Floer homology, Kronheimer-Mrowka's monopole Floer homology is defined for
all closed, oriented 3-manifolds. (The instanton Floer homology has only
been successfully defined for rational homology spheres, except some special cases. However, unlike the closed
4-manifold Seiberg-Witten invariant, which is in most cases
independent of the perturbation form \(\varpi\), the Seiberg-Witten
Floer homology depends on the cohomology class of the perturbation
form in the Seiberg-Witten equation: Let \((M,
\mathfrak{s})\) be a closed spin-c 3-manifold and 
\(X=\mathbb{R}\times M\), and take the form \(\varpi\) in
(\ref{eq:SW}) to be \(\mu+ds\wedge*_3 \mu\) for a closed 2-form on
\(M\), where \(s\) denotes the affine coordinate on the \(\mathbb{R}\)
factor of \(X=\mathbb{R}\times M\), and \(*_3\) denotes the
3-dimensional Hodge dual. The spin-c structure \(\mathfrak{s}\) on \(M\) determines a spin-c
structure on \(X\), which we denote by the same notation. Denote the
associated spinor bundle on \(M\) and \(X\) respectively by
\(\mathbb{S}_X=\mathbb{S}_X^+\oplus \mathbb{S}_X^-\) and
\(\mathbb{S}\). 
Clifford multiplication by \(ds\) determines an isomorphism
\(\mathbb{S}_X^+\simeq \mathbb{S}_X^-\simeq\mathbb{S}\). The equation (\ref{eq:SW}) can be interpreted as a formal
gradient flow equation on \(\op{Conn}(M)\times \Gamma (\mathbb{S})\),
where \(\op{Conn}(M)\) denotes the space of connections on
\(\det \mathbb{S}\):
\begin{equation}\label{eq:SW_s}\frac{d}{ds} (B, \Phi)=-\big(*_3 (F_B-\Phi^\dag \tau\Phi+i\mu), D_B\Phi\big).
\end{equation}
Here, \((B(s), \Phi(s))\) denotes a path in \(\op{Conn}(M)\times
\Gamma (\mathbb{S})\) parametrized by \(s\in \mathbb{R}\). 
The gauge group \(\mathcal{G}:=C^\infty(M, U(1))\) acts on
\(\op{Conn}(M)\times \Gamma (\mathbb{S})\) by \(u\cdot (B,
\Phi)=(B-2u^{-1} du, u\Phi)\) \(\forall u\in \mathcal{G}\) and \((B,
\Phi)\in \op{Conn}(M)\times \Gamma (\mathbb{S})\). 
The  Seiberg-Witten equation (\ref{eq:SW_s}) is invariant under the gauge action,
and thus has an interpretation as the (formal) flow equation of the dual
vector field to a 
closed 1-form on \(\mathcal{B}:=\big(\op{Conn}(M)\times \Gamma
(\mathbb{S})\big)/\mathcal{G}\). Note that
\(\pi_1(\mathcal{B})=H_1(\mathcal{B})=H^1(M;\mathbb{Z})\), and the
cohomology class of the aforementioned closed 1-form is given by
\begin{equation}\label{pert-class}
2\pi^2 c_1(\mathfrak{s})-\pi[\mu]\in
H^2(M;\mathbb{R})\stackrel{P.D.}{\simeq }\Hom
(H^1(M;\mathbb{Z});\mathbb{R})\simeq
H^1(\mathcal{B};\mathbb{R}).
\end{equation}
The Seiberg-Witten-Floer homology associated to
(\ref{eq:SW_s}) is modelled on the Morse-Novikov theory associated this closed
1-form on \(\mathcal{B}\), and thus for this heuristic reason is
expected to  depend on the classes \([\mu]\) 
and \(c_1(\mathfrak{s})\) but not other parameters. This is indeed
verified by a lengthy argument in \cite{KM}. With \(\mathfrak{s}\)
fixed, the class \(-\pi[\mu]\) is called the ``period class'' of the
Seiberg-Witten-Floer homology built from (\ref{eq:SW_s}). 
Meanwhile, it follows from the
Atiyah-Patodi-Singer theorem on spectral flows that this Floer
homology group has a relative grading by the cyclic abelian group
\(\mathbb{Z}/c_{\mathfrak s}\), where \(c_{\mathfrak s}\in
2\mathbb{Z}\) is the gcd of the values of \(c_{\mathfrak s}\), viewed as
a linear map from \(H^1(M;\mathbb{Z})\stackrel{P.D.}{\simeq }
H_2(M;\mathbb{Z})\) to \(\mathbb{Z}\). The grading group is
\(\mathbb{Z}\) when \(c_1(\mathfrak{s})\) is torsion, namely when
\(c_{\mathfrak s}=0\). For technical reasons that we shall not explain
here, given an arbitrary pair of \(c_1(\mathfrak{s})\) and periodic
class \(c\), 
the corresponding Seiberg-Witten Floer homology is only defined for
coefficient rings  \(\Lambda\) that satisfy certain completeness
conditions, called ``\(c\)-completeness'' in \cite{KM}. This Floer
homology is denoted by \(\mathring{HM}(M, \mathfrak{s},
c;\Lambda)\) in \cite{KM}. The coefficient ring \(\Lambda\) is assumed
to be \(\mathbb{Z}\) when \(\Lambda\) is omitted
from the notation, and the period class \(c\) is assumed to be \(0\)
when it is omitted. 
\paragraph{\it Remark.}
To be more precise, for \(HM\) to be well-defined, due to
transversality issues one should allow
certain additional abstract perturbations to
(\ref{eq:SW_s}). Cf. \cite{KM} Chapter 10. We ignore this technical
issue in this article. 
\medskip

\(ECH\) is constructed as a Floer-homology extension of Taubes's
Gromov invariant \(Gr\) in his \(SW=Gr\) theorem
\cite{T} for closed 4-manifolds. An equivalence of \(HM\) and \(ECH\),
in parallel to Taube's theorem for closed 4-manifolds is thus
expected from the very beginning of the development of \(ECH\).  Let
\((M, a)\) be a contact 3-manifold with a contact 1-form \(a\), and
choose the metric on \(M\) so that \(a\) is co-closed. In analogy to
the type of perturbations (\ref{eq:-large_r_pert}) used in Taubes's
proof of his \(SW=Gr\) theorem, choose the 2-form \(\mu\) in
(\ref{eq:SW_s}) to be of the form \(\mu=2r*_3a-iF_{B_K}\) and consider
the associated \(HM\) as \(r\to \infty\). Here, \(B_K\) is the 
connection on \(K^{-1}\) induced from the metric, and \(K:=\ker
a\subset TM\) is equipped with the complex structure given by Clifford
multiplication by \(a\). As the \(Gr\)-side analog of \(HM\), the
chain modules of \(ECH\) are generated by certain union of (weighted)
orbits of the Reeb flow of \(a\) (called ``orbit sets''), and the entries of the
differential as a matrix with respect to the basis consisting of orbit
sets are defined by counting (disjoint, weighted) pseudo-holomorphic
curves in \(\mathbb{R}\times M\) asymptotic to the relevant orbit
sets. It turns out that the actual proof of the equivalence of \(HM\)
and \(ECH\) requires essential new ideas in addition to those in
\cite{T}. See Theorem \ref{thm:Tech} for a more precise statement of
this equivalence theorem, and \cite{T:ECH} for full details. 

Ozsvath-Szabo's \(HF\) may be viewed another \(Gr\) counterpart of
\(HM\), albeit in a less straightforward manner. 
Its motivation comes from a variant of the Atiyah conjecture \cite{A},
in addition to Taubes's philosophy of \(SW\)-\(Gr\)
correspondence. Fix a Heegaard decomposition of a closed 3-manifold
\(M\). Let \(f: M\to \mathbb{R}\) be a self-indexing Morse function adapted to this
Heegaard decomposition. By this we mean that \(f\) has unique maximum
and minimum, and \(G\)-pairs of index 2 and index 1 critical points,
where \(G\) is the genus of the Heegaard surface
\(\Sigma=f^{-1}(3/2)\). Let \(\pmb{\alpha }:=\{\alpha _1, \alpha _2,
\ldots, \alpha _G\}\) denote the set of descending cycles from index 2
critical points on  \(\Sigma\), and let  \(\pmb{\beta }:=\{\beta _1, \ldots, \beta
_G\}\) denote the set of ascending cycles from index 1
critical points on  \(\Sigma\). We call the triple \((\Sigma,
\pmb{\alpha }, \pmb{\beta })\) a Heegaard diagram of \(M\). The idea
of the Atiyah conjecture is to relate the Floer homology of \(M\) to a
Lagrangian Floer homology of \((\mathcal{M}, \mathbb{T}_\alpha,
\mathbb{T}_\beta)\), where \(\mathcal{M}\) is a symplectic manifold
typically coming from the moduli space of a suitable dimensional
reduction of the relevant gauge equation on 3-manifolds, and
\(\mathbb{T}_\alpha\) and \(\mathbb{T}_\beta\) are Lagrangian
submanifolds (typically with singularities) defined from moduli spaces
of solution to the gauge equation on \(f^{-1}[3/2,\infty)\) and
\(f^{-1}(-\infty, 3/2]\) respectively. In the setting of Seiberg-Witten
theory, heuristic reasoning from Taubes's philosophy predicts that for
certain large \(r\) perturbation involving \(*df\), the corresponding
triple  \((\mathcal{M}, \mathbb{T}_\alpha,
\mathbb{T}_\beta)\) should be \((\op{Sym}^G(\Sigma),
\alpha_1\times\cdots\alpha _G, \beta_1\times \cdots \times \beta
_G)\). The Heegaard Floer homology is a variant of Lagrangin Floer
homology associated to  \((\op{Sym}^G(\Sigma),
\alpha_1\times\cdots\alpha _G, \beta_1\times \cdots \times \beta
_G)\), with one extra key ingredient: a
choice of a base point \(z_0\in \Sigma-\bigcup_i\alpha
_i\cup\bigcup_j\beta \). This is used to define a filtration on the
relevant 
Heegaard Floer complex. It is somewhat long and complicated to explain
the relation of \(HF\) with \(HM\) with Taubes's type of
perturbations; the interested reader is referred to \cite{L}.

\section{Equivalences of Floer homologies: theorems and questions}

We may now state the isomorphism theorem of \cite{KLT} more
precisely: 
\begin{theorem}[(\cite{KLT} Theorem V.1.4)]\label{thm:KLT}
Let \(M\) be a closed, oriented 3-manifold, and \(\mathfrak{s}\) be a
spin-c structure on \(M\). Then there exists a system of
isomorphisms from \(HF^\circ_*(M,\mathfrak{s})\), \(\circ=-, \infty, +,
\wedge\), respectively to \(HM^\circ_*(M,\mathfrak{s}, c_b)\),
\(\circ=\wedge, -, \vee, \sim\),  as \(\mathbb{Z}/c_\mathfrak{s} \mathbb{Z}\)-graded 
\(\mathbb{A}_\dag(M)\)-modules, which is natural with respect to the
fundamental exact sequences of the Heegaard and monopole Floer
homologies.   
\end{theorem}

Here, \(c_b\) stands for a ``balanced perturbation'' in the
terminology of \cite{KM}, and refers to the case when the cohomology
class (\ref{pert-class}) is 0.
Among all periodic classes, the balanced case is strongest in the
following sense: it is one for which
the associated \(HM\) can be defined over \(\mathbb{Z}\), and 
the \(c_b\)-completeness condition required for the coefficient ring
of \(HM\) is vacuous. Thus, \(HM(M,
\mathfrak{s}, c; \Lambda)\) of other local coefficients \(\Lambda\) may be
computed via the universal coefficient theorem, together with results in \cite{KM}'s Chapter
31 relating monopole Floer homologies associated to proportional (\ref{pert-class}). It is also the only class for which
\(\overline{HM}\) is nonvanishing, and consequently by the fundamental
exact sequences for \(HM\), the only class
that the two flavors of \(HM\), \(\widehat{HM}\) and
\(\widecheck{HM}\), differ. 

A subtle point worth noting is that in \cite{KM} as well as in other
literature, the monopole Floer homology frequently refers to the ``bullet
version'' (or completed version) \(HM_\bullet\) instead of the ``star
version'' (or pre-completed version) \(HM_*\) appearing in the statement of the
preceding Theorem. The former version, \(HM_\bullet\), uses coefficients that are completed with
respect to the \(U\)-action, and therefore is slightly weaker than the
latter version, 
\(HM_*\).  For example,  \(\overline{HM}_\bullet\) vanishes while
\(\overline{HM}_*\) is nontrivial in the example computed in \cite{KM}
Equation (35.4). Working with \(HM_\bullet\) is in particular more
convenient in discussions involving maps between monopole Floer
homologies induced from cobordisms between two 3-manifolds. However,
to be able to define the Floer chain complex with {\em polynomial} (in \(U\))
coefficients (instead of power series coefficients), certain strong compactness results are necessary.
For \(HM\) and \(HF\), these are guarenateed respectively via the
balanced perturbation condition and the ``strong admissibility'' assumption on the
Heegaard diagram. 

In comparison, the aforementioned finiteness/compactness results are missing in
\(ECH\). This nevertheless is consistent with Taubes's \(HM=ECH\)
theorem, because according to \cite{KM}, \(\widecheck{HM}(M,
\mathfrak{s})\simeq\widecheck{HM}_*(M, \mathfrak{s},
c_b)=\widecheck{HM}_\bullet(M, \mathfrak{s}, c_b)\). Here, \(\mathring{HM}(M, \mathfrak{s})\) stands for the
version of Seiberg-Witten-Floer homology when the perturbation form is
exact. In this case, except for the case when \(c_1(\mathfrak{s})=0\), \(\widehat{HM}(M, \mathfrak{s})=\widecheck{HM}(M,
\mathfrak{s})\), and both of them are finite rank \(\mathbb{Z}\)-modules. 

Meanwhile, just as (the original) \(ECH\) is an analog of
\(\widecheck{HM}\simeq HF^+\), one may define another flavor of
\(ECH\), called 
\(\widehat{ECH}\), as an analog of
\(\widetilde{HM}\simeq\widehat{HF}\). (See e.g. \cite{CGHH}). The pair
\(ECH\) and \(\widehat{ECH}\) fits in an analog of the second long exact sequence in
Equation (\ref{eq:fundsq}) above, also called the fundamental sequence
of ECH. 
It follows directly from Taubes's proof of \(HM=ECH\), together with the definitions of \(\widehat{ECH}(-M, [\xi,
\mathfrak{s}])\) in \cite{CGHH} and \(\widetilde{HM}(M,
\mathfrak{s})\) in \cite{KLT} part V, that the latter is equivalent to
\(\widetilde{HM}\).  

To summarize in a more precise fashion:
\begin{theorem}
\label{thm:Tech}
Let \(M, \mathfrak{s}\) be as in the previous theorem, and let \(\xi\)
be a contact structure on \(M\). Then there is a pair of
isomorphisms \(ECH (-M, [\xi,
\mathfrak{s}])\simeq \widecheck{HM}(M, \mathfrak{s})\); \(\widehat{ECH}(-M, [\xi,
\mathfrak{s}])\simeq \widetilde{HM}(M, \mathfrak{s})\), which  are natural with respect to the fundamental sequences on
both sides. I
\end{theorem}
As a consequence of Theorems \ref{thm:KLT} and \ref{thm:Tech}, one has  
\begin{theorem}\label{thm:CGH}
Let \(M, \mathfrak{s}\) and \(\xi \) be as in the previous
theorem. Then there is a system of isomorphisms  \(ECH (-M, [\xi,
\mathfrak{s}])\simeq HF^+(M, \mathfrak{s})\); \(\widehat{ECH}(-M, [\xi,
\mathfrak{s}])\simeq \widehat{HF}(M, \mathfrak{s})\), which are natural with respect to the fundamental exact
sequences on both sides. 
\end{theorem}
Alternatively, \cite{CGH} has a
``purely symplectic'' proof of the preceding theorem without going through Seiberg-Witten theory.

Comparing Theorems \ref{thm:KLT}, \ref{thm:Tech} and \ref{thm:CGH} above, one
naturally asks: 
\begin{question}
Is there an \(ECH\) analog of (pre-completed) \(\widehat{HM}_*\) or
\(HF^-\), which is defined from Floer chain complexes of
\(\mathbb{Z}[U]\)-modules?
\end{question}

If such an analog can be shown to be equivalent to \(\widehat{HM}\) or
\(HF^-\), the aforementioned finiteness properties of the latter Floer
homologies might help answer questions related to
finiteness of certain types of Reeb orbits on
contact 3-manifolds. 

In a different direction, Taubes's \(HM=ECH\) theorem has a sister version for 3-dimensional
mapping tori, where the contact form is replaced by a  harmonic,
nowhere-vanishing 1-form. A variant of \(ECH\), dubbed ``\(PFH\)''
(periodic Floer homology) by
Hutchings, is shown to be equivalent to (a different version) of
\(HM\) by the present author and Taubes: 

Let \((F, w_F)\) denote a closed oriented surface \(F\) equipped with
a volume form \(w_F\), and \(\varphi\) is a volume
preserving automorphism of \(F\). Let \(M_\varphi\) denote the mapping
torus of \(\varphi\) and \(w_\varphi\) the closed 2-form on
\(M_\varphi\) induced from \(w_F\). Let \(K^{-1}\subset TM_\varphi\) denote the
subbudle consisting of tangent vectors to the fibers of the bundle \(M_\varphi\to S^1\), and \(c_1(K^{-1})\) its Euler class. An element in
\(\Gamma\subset H_1(M_\varphi;\mathbb{Z})\) is said to be {\em monotone}
when \([w_\varphi]=-\lambda (c_1(K^{-1})+2P.D.(\Gamma))\) for a real
number \(\lambda\). It is said to be positive monotone when
\(\lambda>0\), and negative monotone with \(\lambda <0\). The periodic
Floer homology of \(M_\varphi\) in the class \(\Gamma\) is denoted by
\(HP_*(\varphi:(F, w_F)\circlearrowleft , \Gamma)\) in \cite{LT}. See
\S 1.1 therein for details of the definition. In \S 1.2 of the same paper, a
spin-c structure \(\mathfrak{s}_\Gamma\) and a closed 2-form
\(\varpi_r:=2rw_\varphi+\wp\) is assigned to each pair
\((\varphi:(F, w_F)\circlearrowleft , \Gamma)\). Here, \(r\) is a
sufficiently large real number, and \(\wp\) is a closed 2-form in the
cohomology class \(2\pi c_1(\mathfrak{s}_\Gamma)\). The precise choice
of \(r\) and \(\wp\) turns out to be immaterial. 

\begin{theorem}[(\cite{LT}, Theorem 1.1)]\label{thm:PFH}
Let \((F, w_F)\), \(\varphi\) be as above, and let \(\Gamma\) denote
either a positive or negative monotone class. 
Then \(HP(\varphi: (F, w_F)\circlearrowleft  , \Gamma)\simeq HM
(M_\varphi, \mathfrak{s}_\Gamma, [\varpi_r])\). 

 \end{theorem}
The \(HM\) on the right hand side of the isomorphism above stands for
either \(\widehat{HM}\) or \(\widecheck{HM}\), which are the same
under this particular setting. According to \(\cite{KM}\), for a
monotone \(\Gamma\), the periodic class \([\varpi_r]\) is what is
called 
``positive/negative monotone'' with respect to the spin-c structure
\(\mathfrak{s}_\Gamma\) for all sufficiently large \(r\), precisely
when \(\Gamma\) is positive/monotone. In other words, when \(\Gamma\)
is monotone, 
\[
\begin{split}
HP(\varphi: (F, w_F)\circlearrowleft  , \Gamma)
& \simeq HM(M_\varphi,\mathfrak{s}_\Gamma) \quad \text{if \(d_\Gamma:=\langle P.D.[F],
\Gamma\rangle<g-1\), and }\\
HP(\varphi: (F, w_F)\circlearrowleft  , \Gamma)& \simeq HM
(M_\varphi, \mathfrak{s}_\Gamma, c_-) \quad \text{if \(d_\Gamma >g-1\).}
\end{split}
\]

In the above, \(c_\pm\) respectively denotes a positive/negative
periodic class for \(HM\). According to \cite{KM},
\(HM(M_\varphi,\mathfrak{s}_\Gamma, c_+)\simeq
HM(M_\varphi,\mathfrak{s}_\Gamma)\) and \(HM(M_\varphi, \mathfrak{s}_\Gamma, c_-)\) is related to
\(HM(M_\varphi,\mathfrak{s}_\Gamma)\) and
\(\overline{HM}(M_\varphi,\mathfrak{s}_\Gamma, c_b)\) via a long
exact sequence, and therefore is always different from the version of monopole
Floer homology appearing in the \(HM=ECH\) theorem, Theorem \ref{thm:Tech}. (See discussions following Corollaries 1.4 and 1.5 in
\cite{LT}, as well as references given therein). 
 
The contact structure appearing in Theorem \ref{thm:Tech} and the
mapping tori structure appearing in Theorem \ref{thm:PFH} are both
special cases of the so-called ``stable Hamiltonian structure'' (see
e.g. \cite{HT:SHS}). It is therefore natural to ask:
\begin{question}
Is there a version of ECH for 3-manifolds with stable Hamiltonian
structure, that emcompasses both the \(PFH\) for mapping tori and the
\(ECH\) for contact 3-manifolds as special cases?  
\end{question}
As hinted in \cite{CFP}, this might require more than the analytic techniques in \cite{HT}.

Assuming that the preceding question has been positively answered,
\begin{question}
Prove a version of \(HM=ECH\) for this generalized \(ECH\)
that encompasses \(HM=ECH\) of \cite{T:ECH} and \(HM=PFH\) of
\cite{LT} as
special cases.
\end{question}

This is likely a very difficult question that requires essential new ideas. In spite of many
similarities in the proofs of the sister theorems \(HM=ECH\) and
\(HM=PFH\) in \cite{T:ECH} and \cite{LT}, one major difference is the
``energy bound'' (roughly, a certain \(L^1\)-bound of curvature) that
forms the starting point of the proof for geometric convergence of
Seiberg-Witten solutions giving rise to holomorphic curves. In the case of
\(HM=PFH\), it follows mainly from topological reasons; while in the
case of \(HM=ECH\), it follows from a more delicate spectral flow
estimate. Unfortunately, there is no easy way of combining these two very
different arguments in general. 

Some partial results towards this
direction appear in \cite{KLT}. (See in particular papers IV and V of
\cite{KLT}).
As will be explained in Section 3 below in more detail, this series of
articles defines a variant of \(ECH\) for certain manifolds with
stable Hamiltonian structure, and aspects of its relation to
\(HM\) are established. The proofs of these results are long and hard,
and in fact constitute the technical core of the proof of Theorem
\ref{thm:KLT}. However, the method therein work only for a very special
type of stable Hamiltonian structure. (More will be said about this
stable Hamiltonian structure in Section 4 below). Brute force is used to
amalgamate the key energy bounds in the proofs of \(HM=ECH\)
and \(HM=PFH\), that are obtained via very different methods. In these articles, the relevant stable
Hamiltonian structure ``splits'' along certain simple surfaces, such
that one side has a mapping torus structure, and the other a specific
contact structure studied extensively by Taubes. (See
e.g. \cite{T:S2}). The 3-manifold is then stretched very long along
the splitting surface so that the estimates on both sides may be
performed essentially separately.

\section{Local coefficients and general perturbation classes}

No additional work is required to define local coefficients
versions for the Floer homologies \(HM\), \(HF\), \(ECH\). In fact,
the compactness results needed to define a Floer homology are
typically weaker for local coefficients with suitable completeness
properties. In Heegaard Floer theory, this manifects itself in that
only ``weakly admissible'' Heegaard diagrams are required to define 
Heegaard Floer complexes with twisted coefficients. 
 
The isomorphism theorems, Theorem \ref{thm:KLT}, \ref{thm:Tech} and hence
\ref{thm:CGH} extend to isomorphisms in any corresponding local
coefficients, once it is verified that the chain  maps inducing the
isomorphisms in Floer homologies intertwine with the
\(\mathbb{A}_\dag\)-action. E.g., this is done in paper V of \cite{T:ECH}
for the isomorphism \(HM=ECH\). (Very roughly, this is because the
\(H_1(M)/\op{Tor}\) part of the \(\mathbb{A}_\dag\) acts reflects the
action of the relevant ``fundamental group'' on Floer complexes. For an explanation in more
precise terms, see e.g. Section 6 of \cite{LT}). 

A little more needs to be said about defining and extending the
aforementioned isomorphism theorems to general ``perturbation
classes''. (What is what was called the ``periodic class'' in
Seiberg-Witten theory). 

On the \(ECH\) side, a contact form is exact, and the ``perturbation
class'' is always trivial  if one restricts to contact
3-manifolds. ``Perturbed'' versions of ECH enter the scene in 
the more general realm of stable Hamiltonian structures. In
particular, \(PFH\) can be defined for arbitrary \(\Gamma\) and a range of {\em necessarily
nontrivial} perturbation classes and local coefficients that satisfy 
certain completeness conditions. The isomorphism of these with the
corresponding Seiberg-Witten-Floer homology, extending Theorem
\ref{thm:PFH} above, is stated precisely as Theorem 6.5 in
\cite{LT}. The generalized ECH relevant to the proof of Theorem
\ref{thm:KLT} in \cite{KLT} is another example of ECH with nontrivial perturbation class.

On the Heegaard Floer side, a Heegaard Floer homology \(HF^\circ(M,
\mathfrak{s}, \eta)\) corresponding to
other perturbation class \(\eta\in H^2(M, \mathbb{R})\) is defined in \cite{OS} \S 11.0.1, see
also e.g. \cite{Wu} for more details. This ``perturbed Heegaard Floer
homology'' has as coefficients the Novikov ring \(\Lambda_A\)
described presently. This
ring  consists of  formal power series \(\sum_{r\in \mathbb{R}}
a_rT^r\) with \(a_r\in \mathbb{R}\) and \(T\) a formal variable,
satisfying the condition that for any fixed \(N\in \mathbb{R}\), the
number of \(a_r\neq 0\), \(r<N\) is finite. Its
multiplication law is given by \((\sum_{r\in \mathbb{R}}a_rT^r)\cdot (\sum_{s\in \mathbb{R}}
b_sT^s)=\sum_{r\in \mathbb{R}}\sum_{k\in \mathbb{R}}a^kb^{r-k} T^r\).
As a pay off for working with this complicated completed coefficient
ring, no admissibility condition on the Heegaard diagram is required
for the associated Heegaard Floer complex to be well defined. The
cohomology \(\eta\in H^2(M, \mathbb{R})\) enters the definition of the
boundary map \(\partial\) of this perturbed Heegaard Floer complex via
the area assigned to each domain in the Heegaard diagram counted in
\(\partial\).

It takes only superficial
changes (modification of Lemma 1.1 and Lemma 1.2 in paper II) to the proof of Theorem \ref{thm:KLT} in \cite{KLT} to show that:

\begin{theorem}
Let \((M, \mathfrak{s})\) be as in Theorem \ref{thm:KLT}.
There is a system of  isomorphisms from \(HF^\circ_*(M,\mathfrak{s}, \eta)\), \(\circ=-, \infty, +,
\wedge\), respectively to \(HM^\circ_*(M,\mathfrak{s}, c_\eta; \Lambda_A)\),
\(\circ=\wedge, -, \vee, \sim\),  as \(\mathbb{Z}/c_\mathfrak{s} \mathbb{Z}\)-graded 
\(\mathbb{A}_\dag(M)\)-modules, which is natural with respect to the
fundamental exact sequences of the Heegaard and monopole Floer
homologies. 
\end{theorem}

\section{Applications of the equivalence theorems}

The equivalence theorems in Section 1.1 are useful because of the very
different geometric origins of the three Floer homologies. 

Among the aforementioned equivalence theorems, Theorem \ref{thm:Tech}
has more success in finding applications so far. Most notably, the 3-dimensional
Weinstein conjecture and its analogs. See e.g. \cite{T:Weinstein},
\cite{T:Weinstein2}, \cite{HT2}, \cite{HT:SHS}. In another direction,
Hutchings defined an ECH version of symplectic capacities, which
provides complete obstructions to symplectically embedding
4-dimensional ellipsoids to each other. See e.g. \cite{Hs}. These
results make use of the fact that expected properties or definitions in
ECH are often difficult to carry out directly, and thus its relation
with \(HM\) enables one to appeal to the more fully developed
Seiberg-Witten theory. For example, The proof of Weinstein conjecture type
results indirectly make use of the easy computatin of
\(\overline{HM}\); the proof of the Arnold's chord conjecture makes
use of the surgery exact sequences in Seiberg-Witten-Floer homology.
The definition and key properties of ECH capacities go through the
definitions and properties of maps induced by cobordisms in
Seiberg-Witten-Floer theory.  

As immediate consequences of Theorem \ref{thm:KLT}, the equivalence
\(HM=HF\) can be used to compute
one Floer homology of specific 3-manifolds in terms of the other,
depending on which is simpler. For example, the intricate computations of \(HF\)
for mapping tori done in \cite{JM, JM2} and the computation of
\(HF^\infty\) in \cite{M} follow easily from the corresponding
computation on the Seiberg-Witten side. The Seiberg-Witten-Floer
homology of Seibert-fibered spaces has been carried out for general
Seifert fibered spaces in \cite{MOY}, while the compuation of \(HF\)
are only done for various special cases of Seifert fibered spaces in the
literature, using different methods. On the other hand, it has been shown
that \(\widehat{HF}\) is purely combinatorial (see e.g. \cite{SW}),
and thus the equivalence \(HM=HF\) implies that \(\widetilde{HM}\) can
likewise be computed purely combinatorially.

To find more interesting application of  Theorem \ref{thm:KLT} will
most likely require extending the equivalence of {\em Floer
  homologies} of 3-manifolds 
both sides to other aspects of Seiberg-Witten and Ozsvath-Szabo
theories. For example, like the aforementioned applications of
\(HM=ECH\), it is important to know that whether the system of isomorphisms constructed
in \cite{KLT} is natural with respect to the surgery exact sequences
and TQFT structures on both sides.
In the rest of this article, we discuss some possible amplifications
of the arguments in \cite{KLT} towards this goal. These are
emphatically not straightfoward, in particular because the approach
adopted in \cite{KLT} is rather indirect. In preparation, the next
section gives a 
brief outline of the proof in \cite{KLT}. For a fuller summary, see
paper I of \cite{KLT}.

\section{Outlinig the proof of Theorem \ref{thm:KLT}: motivation and strategy}

Surprisingly, the formally similar algebraic structures on
\(HF\) and \(HM\) have completely different origins. For example, The
mechanism that gives rise to the four flavors \(\circ=-, \infty, +,
\wedge\) in \(HF^\circ\), stem from a filtration on the (simplest
flavor of) the Heegaard Floer complex. Roughly speaking, the Floer
complex \(CF^\infty\) can be viewed as a chain complex with local
coefficients \(\mathbb{Z}[\mathbb{Z}]\simeq  \mathbb{Z}[U,
U^{-1}]\). Alternatively, viewing it as an analog of the chain complex
of a \(\mathbb{Z}\)-covering of a finite-dimensional space, the
\(U\)-action corresponds to deck transformations.  A key
consequence of the geometric setup in constructing \(HF\) is that a
basis of the chain module \(CF^\infty\) may be found so that with
respect to which, the boundary map \(\partial^\infty\) has the form of
a matrix with coefficients in the {\em polynomial} ring
\(\mathbb{Z}[U]\). This equips \(CF^\infty\) with a {\em filtration}
by \(\mathbb{Z}[U]\)-subcomplexes \[\cdots CF^-\subset U^{-1}\cdot
CF^-\subset U^{-2}\cdot CF^-\subset\cdots \subset CF^\infty\] 
The four flavors of \(HF^\circ\) are defined as various homology
groups associated to this filtered chain complex. For example, the
first fundamental sequence is the relative long exact sequence induced from the short exact
sequence:
\[
0\to CF^-\to CF^\infty\to CF^+:=CF^\infty/CF^-\to 0.
\]

In contrast, the Seiberg-Witten Floer homology is modelled on an
\(S^1\)-equivariant homology theory. From this point of view, the
module structure over \(\mathbb{Z}[U]\) in \(HM\) reflects the
module-structure of \(S^1\)-equivariant homologies. (\(\mathbb{Z}[U]\)
is the cohomological ring of the classifying space of \(S^1\) actions,
\(\mathbb{C}P^\infty\)). The four flavors of \(HM\)  arise as various
versions of \(S^1\)-equivariant homologies. For example, the first
fundamental sequence is modelled on a well-known long exact sequence in
\(S^1\)-equivariant homology theory relating the Borel equivariant
homology, the Tate version of equivariant
homology, and the so-called co-Borel version.  

To resolve the seemingly irreconcilable differences in the foundations
of Heegaard and Seiberg-Witten Floer homologies,
the strategy of \cite{KLT} is to go
through a third, intermediate version of Floer homology, denoted
\(ech^\circ\) in paper I of \cite{KLT}. The isomorphisms in
\ref{thm:KLT} will follow from composing  a numbef of isomorphisms
among several Floer homologies.

Consider \(\mathring{HM}(M, \mathfrak{s}, c_b)\) and \(HF^\infty(M,
\mathfrak{s})\)  from the statement of Theorem \ref{thm:KLT}. The
relevant intermediate Floer homology \(ech^\circ\) is defined for an
auxiliary 3-manifold \(Y:=Y_M\). Given a pointed Heegaard diagram used
to define \(HF^\infty(M, \mathfrak{s}), \) let \(f\) be a
self-indexing Morse function associated with this Heegaard
decomposition. Suppose this \(f\) has one pair of index 0 and index 3
critical points, and \(G\) pairs of index 1 and index 2 critical
points. \(Y_M\) is built from \(M\) and by doing a 0-dimensional
surgery along the aforementioned \(G+1\) copies of \(S^0\)'s
(i.e.pairs of critical points). Denote the copy of \(I\times
S^2\subset Y_M\) resulted from surgery along the pair of index 0 and
index 3 critical points by \(\mathcal{H}_0\), and denote by
\(\mathcal{H}_i\), \(i\in \Lambda:=\{1, \ldots, G\}\), each of the \(G\) copies of \(I\times
S^2\subset Y_M\) resulted from surgery along a pair of index 2 and
index 1 critical points. 
This \(Y_M\) is then assigned with a special stable Hamiltonian
structure,  \(a\in \Omega^1(Y_M)\) and \(w\in \Omega^2(Y_M)\), such that \(a\wedge w\) is nowhere vanishing,
\(dw=0\), and \(da=\lambda w\) for a scalar function \(\lambda\) on
\(Y_M\), \(\lambda\geq 0\). Choose a metric on \(Y_M\) such that \(a\)
agrees with the Hodge star of \(w\). The salient features of this
special stable Hamiltonian structure  includes: 
\begin{itemize}
\item The metric is such that the interval factor \(I\) in
  \(\mathcal{H}_i\simeq I\times S^2\) is very long cylinder.
\item \(\lambda\) ``approximates'' the zero function in the interior
  of \(Y_M-\bigcup_{i\in \Lambda}\mathcal{H}_i\). In other words, both
  \(a\) and \(w\) are ``almost harmonic'' in this region.
\item \(\lambda\) ``approximates'' a nonzero constant function in the interior
  of each \(\mathcal{H}_i\), \(i\in \Lambda\). In other words, the
  1-form  \(a\) is ``almost contact'' on each \(\mathcal{H}_i\).
\item Over the long cylinder \(\mathcal{H}_0\), \(a\) approximates the
  harmonic form \(dt\) on \(\mathbb{R}\times S^2\), where \(t\) is an
  affine parameter of the first factor \(\mathbb{R}\).
\item Over the interior of \(Y_M-\bigcup_{i\in
    \{0\}\cup\Lambda}\mathcal{H}_i\), \(a\) approximates a constant
  multiple of \(df\).
\item Over each long cylinder \(\mathcal{H}_i\simeq I\times S^2\),
  \(i\in \Lambda\), \(a\) approximates the following well-known
  contact form \(a_0\) on \(\mathbb{R}\times  S^2\):
  \[a_0=(1-3\cos^2\theta)dt-\sqrt{6}\cos\theta\sin\theta^2 d\phi, \]
where  \((t, \theta, \phi)\) is the standard cylindrical coordinates of
  \(\mathbb{R}\times S^2\).
\end{itemize}

The contact form \(a_0\) above has an analog on the closed 3-manifold
\(S^1\times S^2\), in which case \(t\in S^1\) instead of
\(\mathbb{R}\). The p-holomorphic curves in the symplectization of
this contact 3-manifold have previously been studied extensively by Taubes, see
e.g. \cite{T:S2}. The contact manifold \(S^1\times S^2\) is of
particular interest because its symplectization serves as an
asymptotic model for the (real blow-up of) tubular neighborhoods of
connected 
components of zero loci of generic self-dual harmonic 2-forms on
closed 4-manifolds. A harmonic 2-form of this type is regarded as a generalization
of symplectic forms; it is sometimes called a ``singular symplectic
form''. Unlike symplectic forms, it exists on any
closed 4-manifolds with \(b_2^+>0\). See e.g. \cite{Ts} for an
explanation of the relevance of such harmonic
2-forms in 4-manifold topology.

Roughly speaking, singular symplectic forms arise in the setup of \cite{KLT}
in the following way: Let \(\underline{M}\) denote the 3-manifold
obtained from \(M\) by performing a 0-dimensional surgery along the
pair of index 0 and index 3 critical points of \(f\). The real-valued
Morse function \(f\) extends naturally to an \(S^1\)-valued Morse
function \(\underline{f}\) on \(\underline{M}\), which has no
extrema. A result of Calabi asserts that a metric on \(\underline{M}\)
may be found with respect to which \(\underline{f}\) become
harmonic. This naturally gives rise to a singular symplectic form on
the cylinder \(\mathbb{R}\times \underline{M}\) with product metric,
namely \[\omega=ds\wedge d\underline{f}+*_3d\underline{f},\]
where \(s\) is an affine coordinate for the \(\mathbb{R}\)-factor of
the cylinder \(\mathbb{R}\times \underline{M}\), and \(*_3\) denotes
the 3-dimensional Hodge star. In view of the above discussion, the
construction of \(Y_M\) and the special Hamiltonian structure \((a,
w)\) should be understood as doing real blow-ups of
\(M\) along critical points of \(f\), then ``connect pairs of ends'' of this
blown-up manifold to get a closed manifold. There is a good reason for
choosing to work with closed 3-manifolds instead of manifold with cylindrical ends:
The analysis required for constructing Floer homologies of
noncompact manifolds is almost always intractible.

Fix a spin-c structure \(\mathfrak{s}\) on \(M\). 
There is a natural way to choose a corresponding spin-c
structure on \(Y_M\). As explained e.g. in \cite{HL2}, one
may represent \(\mathfrak{s}\) by a set of mutually disjoint embedded arcs
\(\{\gamma_i\}_{i\in \{0\}\cup \Lambda}\), where \(\partial\gamma_i\simeq S^0\) is the attaching sphere for
\(\mathcal{H}_i\). 
Meanwhile, by ``pinching'' along the boundary
2-sphere of a tubular neighborhood of each \(\gamma_i\), \(Y_M\) may
be expressed as a connected sum of \(M\) with \(G+1\) copies of
\(S^1\times S^2\), one for each \(\mathcal{H}_i\). The spin-c
structure on \(Y_M\) corresponding to \(\mathfrak{s}\) is, with
respect to this connected sum decomposition, the connected sum of of
\(\mathfrak{s}\) for the \(M\)-summand, together with
\(\mathfrak{s}_0\) for each copy of \(S^1\times  S^2\) corresponding
to \(i\in \Lambda\), and with \(\mathfrak{s}_K\) on the copy of
\(S^1\times S^2\) corresponding to \(i=0\). Here, \(\mathfrak{s}_0\)
denotes the trivial spin-c structure, namely, the spin-c structure
with trivial first Chern class. The spin-c structure
\(\mathfrak{s}_K\) is the one represented by the oriented 2-plane field
\(\ker (dt)\), \(t\in S^1\) being an affine parameter of the first factor
of \(S^1\times S^2\).  Fix a choice of the arcs \(\{\gamma_i\}_i\),
and write the corresponding connected sum decomposition of \(Y_M\) as
\begin{equation}\label{eq:Y_as_sum}
Y_M\simeq M\#_{\gamma_0}(S^1\times S ^2)\#_{\gamma_1}(S^1\times S ^2)\cdots\#_{\gamma_G}(S^1\times S ^2).
\end{equation}
With respect to this connected sum decomposition, the graded algebra
\(\mathbb{A}_\dag(Y_M)\) is identified with a tensor product of the
graded algebras
\(\mathbb{A}_\dag(M)\otimes\bigotimes_{i\in
  \{0\}\cup\Lambda}\bigwedge_i\), where each
\(\bigwedge_i=\bigwedge^*\mathbb{Z}_{(-1)}\), the total exterior
algebra on the free graded \(\mathbb{Z}\)-module on a single
generator of degree \(-1\). The latter generator arises from a generator
of \(H_1(S^1\times S^1)/\text{tor}\). The following subalgebra of
\(\mathbb{A}_\dag(Y_M)\) will be useful later:
\(
\hat{\mathbb{A}}_\dag(M)=\mathbb{A}_\dag(M)\otimes\bigotimes_{i\in\Lambda}\bigwedge_i
\)
with respect to the aforementioned factorization of \(\mathbb{A}_\dag(Y_M)\).

It is shown in paper I of \cite{KLT} that one may define a variant of (filtered)
\(ECH\) on \(Y_M\) with the special stable Hamiltonian structure
\((a,w)\). This \(ECH\) comes in four flavors, and is denoted as a
whole as \(ech^\circ\). Like the Heegaard Floer homology, the
superscript \(\circ\) stands for \(-, \infty, +, \wedge\), and they
fit into fundamental long exact sequences in exact analogy with the 
Heegaard Floer homology. They also come equipped with 
\(\mathbb{A}_\dag (Y_M)\)-module structures. The origin of the four
flavors and the \(U\)-action are in complete paprallel to those in
\(HF\).

By taking the cylinders \(\mathcal{H}_i\) in \(Y_M\) to be long and
thin, the decomposition \(Y_M\) as a union of \(M_\delta\) (\(M\) with
small balls near critical points removed) and all the
\(\mathcal{H}_i\)'s, a gluing argument (cf. papers II and III of \cite{KLT})
allows one to compute \(ech^\circ \) by combining arguments from
prior work of Lipschitz and Taubes: In \cite{Li}, Lipshitz reinterpreted
the Heegaard Floer homology as a certian variant of ECH (with Lagrangian
boundary condition) on \(I\times \Sigma\), where \(\Sigma \) is the
Heegaard surface, and the Lagrangian boundary condition given by 
curves in the Heegaard diagram. 
On the other hand, Taubes's work in \cite{T:S2} and its 
sequels give explicit description of holomorphic curves in the
symplectization of \(S^1\times S^2\) equipped with the contact
structure previously mentioned. The result of this computation is
summarized as follows:

Let \(\hat{V}=\mathbb{Z}[y]\) with an odd generator \(y\) of degree
\(1\). Regard this as a free \(\bigwedge^*\mathbb{Z}_{(-1)}\)-module with the
degree \(-1\) generator of the graded algebra
\(\bigwedge^*\mathbb{Z}_{(-1)}\) acting as \(\partial_y\).
 \begin{theorem}
There exists a system of graded \(\hat{\mathbb{A}}_\dag (M)\)-module isomorphisms from \(ech^\circ\) to
  \(HF^\circ(M, \mathfrak{s})\otimes _{\mathbb{Z}}\hat{V}^{\otimes
    G}\). 
These isomorphisms are natural with respect to the fundamental sequences of \(HF^\circ\).
\label{thm:A}
\end{theorem}

This is a direct consequence of Theorem 1.1 in paper III of
\cite{KLT}, and constitutes the first system of isomorphisms for the proof of Theorem
\ref{thm:KLT}. 

The second system of isomorphisms is the stable
Hamiltonian analog of Theorems \ref{thm:Tech} and \ref{thm:PFH}
mentioned in Section 1, and is summarized in the next theorem. There, \(H^\circ_ {SW}(Y)\)
denotes a filtered version of large perturbation Seiberg-Witten Floer
homology on the special stable Hamiltonian 3-manifold \(Y\). 

\begin{theorem}
[(Theorem 3.4, paper I of \cite{KLT})]
There exists a system of isomorphisms from \(ech^\circ\) to
  \(H^\circ_{SW}(Y)\) as \(\mathbb{A}_\dag(Y)\)-modules that is natural with respect to the
  fundamental exact sequences on both sides.
\label{thm:B}\end{theorem}
See paper IV of \cite{KLT} for the definition of \(H^\circ_{SW}(Y)\)
as well as a complete proof. 

The third isomorphism theorem used for the proof of Theorem
\ref{thm:KLT} relates the filtered large perturbation
Seiberg-Witten-Floer homologies \(H^\circ_{SW}(Y)\) with the balanced
Seiberg-Witten-Floer homologies of \(M\). This is summarized in the next
theorem.  The notion \(\op{HMT}^\circ\) below refers a filtered version of
large-perturbation Seiberg-Witten-Floer homology on
\(\underline{M}\simeq M\#_{\gamma_0}(S^1\times S^2)\), originally
introduced in \cite{L}.

\begin{theorem}[(Theorem 1.1 of paper V of \cite{KLT})]\label{thm:conn_sum}
{\bf (1)} There exists a system of 
isomorphisms of \(\hat{\mathbb{A}}_\dag(M)\)-modules 
\[
H_{SW}^\circ(Y)\stackrel{\simeq}{\longrightarrow}\op{HMT}^\circ\otimes\hat{V}^{\otimes G}, \quad \circ=-, \infty, +, \wedge, \]
that preserves the relative
gradings and is natural with respect to the fundamental long exact
sequences on both sides. 

{\bf (2)} There exists a system of canonical isomorphisms of
\(\mathbb{A}_\dag(M)\)-modules from
\[\text{\(\op{HMT}^\circ\), \(\circ=-, \infty, +, \wedge\) respectively to
\(\mathring{HM} \,( M, \mathfrak{s}, c_b)\), \(\circ=\wedge, -, \vee, \sim\),}\]
that preserves the relative
gradings and is natural with respect to the fundamental long exact
sequences on both sides. 
\end{theorem}

An ingredient of the proof  of part (2) above involves some homological
algebra related to the so-called ``Kozsul duality'' in
\(S^1\)-equivariant homology theories. This mechanism converts the
algebraic structures of \(HF^\circ\) and \(HMT^\circ\), resulting from
filtration, to the algebraic structures on \(\mathring{HM}\), resulting
from \(S^1\)-equivariance. 

A central part of  the proof consists of certain filtered
versions of connected sum formula for perturbed Seiberg-Witten-Floer
homologies.  The filtered cobordism formula makes use of cobordism
maps between Floer homologies of two stable Hamiltonian 3-manifolds
\(Y_-\), \(Y_+\) induced from a cobordism \(W\) corresponding to attaching a 4-dimensional
3-handle along a separating 2-sphere in \(Y_-\) and their ``time-reversal'' \(\bar{W}\)
(corresponding to attaching 4-dimensional 1-handles along two points
lying on different connected components of the \(Y_+\)).   The filtration on the perturbed Seiberg-Witten Floer
complex of \(Y_\pm\) is induced by a particular Reeb orbit \(\gamma_\pm\) of the stable
Hamiltonian structure. For example, for \(Y_\pm=Y_M\),
\(\gamma_\pm\) is the curve \(\gamma_0\) through \(\mathcal{H}_0\)
described before. It is essential to show that the aforementioned cobordism maps
perserve the filtration, and this constitutes the technical core and
occupies the bulk of paper V of \cite{KLT}. (In fact, a special case
is required to show that \(H_{SW}^\circ\) and \(HMT^\circ\) are
well-defined filtered Floer homologies, and this occupies a
substantial part of paper IV of \cite{KLT}). To this end, a particular 
holomorphic cylinder \(\mathbb{R}\times S^1\) in \(W\) ending in
\(\gamma_-\) and \(\gamma_+\) is introduced,
and the key is a positivity result of certain curvature
integral over this cylinder. Cf. Proposition 3.4 in paper V of  \cite{KLT}.

It is important to note that the aforementioned positivity result
applies only to very special cobordisms. For example, \(W\) must
contains an open set diffeomorphic to \(\mathbb{R}\times S^2\). A  set of very
stringent conditions on the asymptotic behavior of the geometry
of \(W\) is also required. Cf. Sections 3.2-3.3 of paper V of
\cite{KLT}. As a result, a rather lengthy portion of the latter article is
spent on constructing cobordisms satifying these conditions. (Section
9 of \cite{KLT}, paper V). These special cobordisms will be useful for
comparision the 4-dimensional aspects of the Seiberg-Witten and
Heegaard Floer theories. More about this will be said in the ensuing
sections. 

It is highly desirable to generalize the above positivity result to
more general cobordisms. For general 4-dimensional cobordisms equipped
with a nontrivial self-dual harmonic form, this is carried out in
\cite{L:har}.

\section{4-manifolds and TQFT}

According to Atiyah \cite{A}, 3- and 4-manifold gauge invariants
should fit into a certain ``Topological Quantum Field Theory'' (TQFT),
which assigns (Floer) homology groups to closed 3-manifolds, and for
4-dimensional cobordisms, homomorphisms
between the Floer groups of the 3-manifolds at the beginning
and the end of the cobordism. Aside from some minor glitches (to be
explained later), the TQFT structures on both Heegaard
Floer theory and Seiberg-Witten theory have been fairly well
developed.  In contrast, due to technical
difficulties, there is little progress on constructing TQFT for ECH. 
\begin{question}
Given a symplectic cobordism \((W, \omega)\) between two 3-manifolds
with contact forms \((Y_0, \alpha_0)\), \((Y_1, \alpha_1)\), are there
appropriate ECH maps associated to \((W, \omega)\) (defined by counting
p-holomorphic curves in \((W, \omega)\))? Show that these are equivalent
to the Seiberg-Witten cobordism map \(\widecheck{HM}(W)\). 
\end{question}
Regarding the second question, there are partial results towards the
\(HM\to ECH\) direction of the proposed equivalence for exact
symplectic cobordisms in \cite{HT2}. A parallel question may be asked
for PFH: Consider the category \(\mathcal{C}\) whose objects are pairs \((Y,
\theta)\), where \(\theta \) is a nowhere vanishing harmonic 1-form on
the 3-manifold \(Y\), and whose morphisms between \((Y_0, \theta_0)\),
\((Y_1, \theta_1)\) are pairs \((W, w)\)
consisting of a 4-dimensional cobordism \(W\) and a nowhere vanishing
harmonic 2-form \(w\) on \(W\) satisfying the following conditions:
\((W, w)\) are asymptotic to pairs \((Y_0, w_0)\), \((Y_1, w_1)\) in the
ends, where \(w_0=*\theta_0\), \(w_1=*\theta_1\) are nowhere vanishing harmonic 2-forms on
the 3-manifolds \(Y_0\), \(Y_1\) respectively. By adapting the works
by Donaldson and Gompf relating Lefschetz fibrations and symplectic
structures (cf. e.g. \cite{D}), it should not be difficult to
demonstrate that this category is equivalent to the category called
``FCOB'' in \cite{U}.
\begin{question}
Construct a TQFT from the category \(\mathcal{C}\) that associates to
each object \((Y, \theta)\) its PFH, and to each morphism \((W, w)\) a map between the PFH of the
ends of the cobordism. Show that this TQFT is isomorphic to a
restriction of the (large perturbation) \(HM\) TQFT. 
\end{question}
An application of the convergence theorem in \cite{L} provides the
\(HM\to PFH\) part of the solution to the second question
above. As for the first question, an analog of the desired \(PFH\) TQFT is
constructed by Usher in \cite{U}; in fact, they are expected to be equivalent. 

We now provide a little more details about the 4-manifold
invariants in \(HM\) and \(HF\) theories. 
Let \(Y_0\), \(Y_1\) be {\em nonempty} closed, {\em connected}, oriented 3-manifolds,
and let \(W\) be a connected oriented cobordism from the former to the
latter. Fix a spin-c structure \(\mathfrak{s}\) on \(W\); let
\(\mathfrak{t}_0:=\mathfrak{s}_W |_{Y_0}\),
\(\mathfrak{t}_1:=\mathfrak{s}_W|_{Y_1}\) respectively.  
In \cite{OS:4d}, Ozsvath-Szabo defined 
maps \(F^\circ(W, \mathfrak{s}): HF^\circ(Y_0, \mathfrak{t}_0)\to
HF^\circ (Y_1, \mathfrak{t}_1)\) for \(\circ=-,
\infty, +, \wedge\) along the following lines: First, decompose \(W\)
into a sequence of elementary cobordisms, each corresponding to either
a 1-, 2-, or 3-handle attachment. We call these respectively the index
1, 2, or 3 elementary cobordisms below. Explicit formulae are given for
cobordism maps associated to each type of elementary cobordisms. The cobordism map
\(F^\circ(W)\) is then defined as the composition of cobordism maps
associated with the each step of the aforementioned handlebody
decomposition of \(W\). 

Let \((X, \mathfrak{s})\) be a {\em closed} connected spin-c 4-manifold with \(b_2^+>1\).
Ozsvath-Szabo also introduced an invariant \(\Phi_{X, \mathfrak{s}}\) taking
values in \(\mathbb{Z}/\pm\) from a mixture of the aforementioned
cobordism maps in different flavors: 
Take out two balls in \(X\) and view the latter as a cobordism \((W_X,
\mathfrak{s}_W)\) from \(S^3\) to
\(S^3\). It is shown (cf. \cite{OS:4d}
Lemma 8.2) that \(F^\infty(W, \mathfrak{s})=0\) for any  connected
oriented cobordism \((W, \mathfrak{s})\) between connected oriented
closed 3-manifolds \((Y_0, \mathfrak{t}_0)\),
\((Y_1,\mathfrak{t}_1)\), if \(b^+_2(W)>0\). Thus, in this case the map
\(F^-(W, \mathfrak{s})\) has a lift, denoted \(F^{mix}(W, \mathfrak{s}): HF^-(Y_0, \mathfrak{t}_0)\to
HF^+(Y_1, \mathfrak{t}_1)\). Moreover, when
\(b^+_2(W)>1\), this lift is canonical. Noting that both the Heegaard Floer
homologies \(HF^-(S^3)\) and \(HF^+(S^3)\) are isomorphic to
\(\mathbb{Z}\), denote a generator of the former by \(\Theta_+\),
and 
a generator of the latter by \(\Theta_-\). Set \(\Phi_{X,
  \mathfrak{s}}\) to be the coefficient of \(\Theta_+\) in \(F( W_X, \mathfrak{s}_W)\).

The 4-manifold story on the Seiberg-Witten side is completely parallel.
For each flavor \(\circ=\wedge, -, \vee, \sim\),
\cite{KM} defines a homomorphism \[\mathring{HM} (W, \mathfrak{s}):
\mathring{HM}_\bullet(Y_0, \mathfrak{s}_0)\to
\mathring{HM}_\bullet(Y_1, \mathfrak{s}_1),\]
and for \(W\) with \(b_2^+>1\), there is also a mixed invariant
\(\vec{HM} (W, \mathfrak{s})\). Let \(1\), \(\check{1}\) respectively
denote the standard generators of \(\widehat{HM}_\bullet (S^3)\simeq \mathbb{Z}\)
and \(\widecheck{HM}^\bullet (S^3)\simeq \mathbb{Z}\), it is shown in Proposition 27.4.1 (cf. also Propositions
3.6.1 and 3.8.2) of \cite{KM} that the coefficient of the mixed
invariant, \(\langle \vec{HM} (W_X, \mathfrak{s}_W)1, \check{1})\) is
equal to the closed 4-manifold Seiberg-Witten invariant
\(\mathfrak{m}(X, \mathfrak{s})\) introduced in late 1990's. The
latter is defined by counting solutions of the Seiberg-Witten
equations on \(X\); cf. also \cite{KM} Definitions 27.1.6-1.7. 

With these invariants defined, 
Corollary 23.1.7 in \cite{KM} asserts that the Seiberg-Witten theory
is a weak version of TQFT. It is not a TQFT in the strictest sense, as
instead of the category of all oriented 3-manifolds and 4-dimensional
cobordisms, the Seiberg-Witten functor is defined on the smaller category, denoted
by \(\textsc{COB}\) in \cite{KM}, whose objects consist of {\em
  nonempty, connected} oriented 3-manifolds. In this TQFT the Floer homologies come in
four flavors, and duality takes Floer homologies in one flavor to that
of {\em another}.  As a consequence, the additional condition
\(b_2^+>1\) is required to recover the closed 4-manifold invariant
from the Seiberg-Witten TQFT. The Heegaard Floer theory has a parallel
(weak) TQFT structure, established recently in \cite{JT}.
It is natural to expect: 
\begin{conjecture}\label{conj}
The isomorphisms between \(HF\) and \(HM\) in Theorem \ref{thm:KLT} are 
natural with respect to the TQFT structures on both sides, and that
their closed 4-manifold invariants agree. 
\end{conjecture}
If established, this would
imply the equivalence of their respective contact invariants, and that
the \(HF=HM\) isomorphism is natural with respect to the surgery exact
sequences on both sides. 

Unlike Heegaard Floer theory, the sign ambiguity in the Seiberg-Witten 
4-manifold invariants \(\mathring{HM} (W, \mathfrak{s})\) and
\(\mathfrak{m}(X, \mathfrak{s})\) can be eliminated by fixing
``homological orientations''. This leads one to ask:
\begin{question}
Can the Heegaard Floer invariants \(F_W^\circ\) and
\(\Phi_X\) be refined to obtain \(\mathbb{Z}\)-valued invariants for fixed homological orientations?
\end{question}
There is brief discussion in \cite{KM} generalizing the
Seiberg-Witten 4-manifold invariants above to non-exact perturbations
and local coefficients. It would be interesting to explore the 
corresponding ``perturbed Heegaard Floer 4-manifold invariants''
extending the ``perturbed Heegaard Floer homology'' mentioned in
Section 2 above to a TQFT. 

Meanwhile, one may try to remove the
connectedness assumption in the definition of cobordism maps in either
theory:
\begin{question}
What would be an appropriate general formulation for invariants of 4-manifold
cobordisms between possibly disconnected 3-manifolds, in either
Seiberg-Witten or Heegaard Floer theory?
\end{question}
Cobordism maps between possibly disconnected 3-manifolds has been
defined in some simple 
special cases, such as those used for the connected sum formula
mentioned in Section 4.

\section{Sketching a proof of Conjecture \ref{conj}}

In an article under preparation, the author will give a proof of
Conjecture \ref{conj}. Some salient features of this proof, especially
those pertaining to the prior discussion, are summarized below.

First, in view of the construction of the 4-manifold invariants
  \(F^\circ_{W, \mathfrak{s}}\) and \(\Phi_{X, \mathfrak{s}}\) from
  compositions of elementary cobordisms and the composition theorems
  of their counterparts in Seiberg-Witten theory (cf. \cite{KM}), it suffices to
  compare the cobordism maps on both sides for elementary
  cobordisms. 

It is verified in \cite{KLT} that the isomorphisms in Theorem \ref{thm:KLT}
 map the aforementioned generators \(\Theta_-\in HF^-(S^3)\), \(\Theta_+\in
 HF^+(S^3)\) respectively to the generators \(1\in
 \widehat{HM}(S^3)\), \(\check{1}\in \widecheck{HM} (S^3)\) modulo
 signs. 

Index 1 and index 3 elementary cobordisms are ``time-reversals''
  of each other; so we consider only one of them. Take an index
  1 cobordism \(W\) from \(M\) to \(M'\simeq M\# (S^1\times
  S^2)\). Choose a Heegaard diagram \((\Sigma, \pmb{\alpha}, \pmb{\beta})\) for \(M\) and let \(f\) be a
  self-indexing Morse function on \(M\) associated to this Heegaard
  diagram. The chosen Heegaard diagram for \(M\) induces one for
  \(M'\), \((\Sigma', \pmb{\alpha}', \pmb{\beta}')=(\Sigma,
  \pmb{\alpha}, \pmb{\beta})\#(E, \alpha_0, \beta_0)\), where \((E,
  \alpha_0, \beta_0)\) denotes the so-called standard Heegaard diagram
  of \(S^1\times S^2\), with the
  Heegaard surface \(E\) being a torus, and \(\alpha_0, \beta_0\) are
  embedded circles on \(E\) intersecting transversely at two
  points. Correspondingly, this Heegaard diagram is associated to a
  self-indexing Morse function \(f'\) on \(M'\) that has a single pair
  of index 3 and index 0 critical points, and \(G+1\) pairs of index 2
  and index 1 critical points. Here, the first \(G\) descending cycles
 from index 2 critical points are the \(G\) mutually disjoint circles
  \(\pmb{\alpha}\) on the \(\Sigma\)-summand of \(\Sigma'\), and the first \(G\) ascending cycles from the index 1 critical points are the \(G\) mutually disjoint circles in
  \(\pmb{\beta}\) on the \(\Sigma\)-summand of \(\Sigma'\).  The
  descending and ascending cycles of the last pair of index 2-index 1
  critical points lie on the \(E\)-summand of \(\Sigma'\) and are
  respectively \(\alpha_0\) and \(\beta_0\). By expressing the
  Heegaard diagram of \(M'\) as a connected sum in this way, the
  chain group \(CF^\circ(M')\) may be expressed as a tensor product
  \(CF^\circ(M)\otimes \hat{V}\), and  in \cite{OS:4d}, the map
  \(F^\circ_W: HF^\circ(M)\to HF^\circ (M')\) is defined as the
  map induced by \(CF^\circ(M)\to CF^\circ(M)\otimes \hat{V}:
  \xi\mapsto \xi\otimes y\).  

Recall from Section 4 the construction of the auxiliary manifold \(Y\) with
a stable Hamiltonian structure from a
3-manifold \(M\) and a self-indexing Morse function \(f\) on it. As this
construction will be applied to different pairs of 3-manifolds and
Morse functions, to be specific we denote the stable
Hamiltonian structure constructed from the
pair \((M, f)\) by \(Y(M, f)\). In paper V of \cite{KLT}, as
intermediate steps for the application of filtered connected sum
formula, we also constructed a family of related auxiliary 3-manifolds \(Y_i\), \(i=0,
\ldots, G\), so that \(Y_0=Y\) and \(Y_G=\underline{M}\). These will
be denoted by \(Y_i(M, f)\). Topologically, \(Y_i\) is diffeomorphic
to \(M\) connected summing with \(G+1-i\) copies of \(S^1\times
S^2\). The geometric structures on \(Y_i's\) for \(i>0\) are not
stable Hamiltonian, but like the pair \((a, w)\) of 1- and 2-forms
characterizing a stable Hamiltonian structure, they are associated
with certain pairs of 1- and 2-forms \((a_i, w_i)\) that may vanish
somewhere on \(Y_i\). Each of these \(Y_i\)'s contains a circle used
to define a filtration on the associated large-perturbation
Seiberg-Witten Floer complex. This circle lies away from the zeros of
\((a_i, w_i)\), and abusing notation, we denote it by \(\gamma_0\) for all \(i=0,
\ldots, G\).

Perturb the Morse function \(f\) inside a small ball
\(B^3_{\mathfrak{q}}\) to create a pair
of cancelling index 1 and index 2 critical points \({\mathfrak{q}}\), and denote the
resulting Morse function \(f_s\). 
By a simple modification of the construction in Section 9 of paper V
of \cite{KLT}, we
construct an index 1 elementary cobordism \(W_i\) between \(Y_{i-1}(M,
f)\) and \(Y_{i}(M', f')\), with \(Y_{-1}(M, f):=Y(M, f_s)\). It is 
endowed with a geometric structure so that the associated cobordism
map preserves filtration. 
By composing these index 1 elementary cobordisms with the appropriate
index 1 cobordisms in the proof of the filtered connected
sum formula, and noting that the order of the composition may be
permuted, by induction the computation of \(\mathring{HM}(W)\) may be
reduced to the computation of the filtered cobordism map associated to
\(W_0\). The computation of the latter makes further use of 
composition results: Compose the index 1 elementary
cobordism \(W_0\) with an index 3 elementary cobordism \(W'\) from the proof
of the filtered connected sum formula, with the attaching sphere being
\(\partial B^3_{\mathfrak{q}}\). The right end of this cobordism
is a disjoint union of \(Y=Y(M, f_s)\) with \((S^1\times S^2)\), the
latter being equipped with the round metric and an exact
perturbation. Attach a 2-handle along the core circle in \(S^1\times
S^2\) to form a final cobordism which is diffeomorphic to the
product cobordism  \(\mathbb{R}\times Y\) with a 4-ball
removed. Evaluating the associated cobordism map on the generator
\(\check{1}\) of \(\widecheck{HM}(S^3)\) then yields the identity map
on the Seiberg-Witten-Floer homology of \(Y\). With this, the
cobordism maps associated to \(W_0\) may be computed from those of
\(W'\), which is known from the filtered connected sum formula.

Next, we compare maps associated to index 2 elemenary cobordisms. Let
\(W_K\) now be an elementary cobordisms from \(M\) to \(M_K\), obtained
by attaching a 4-dimensional 2-handle along a framed knot \(K\subset
M\). In \cite{OS:4d}, the cobordism map \(F_{W_K}^\circ\) is defined
by counting ``holomorphic triangles''. Take a Heegaard diagram
\((\Sigma, \pmb{\alpha}, \pmb{\beta})\) adapted to the pair \((M,
K)\). This means that \(M\) is obtained by gluing two
\(G\)-handlebodies along the sets \(\pmb{\alpha}=\{\alpha_1, \ldots,
\alpha_G\}\), \(\pmb{\beta}=\{\beta_1, \ldots, \beta_G\}\) of disjoint
\(G\)-circles, where \(\beta_G=m_K\), the meridian of the knot
\(K\). Meanwhile, the framing of \(K\) is represented by a circle
\(\gamma_G\) on \(\Sigma\). Let \(\pmb{\gamma}=\{\gamma_1, \ldots,
\gamma_G\}\), where for \(i=1, \ldots, G-1\), \(\gamma_i\) is a circle
on \(\Sigma\) obtained by perturbing \(\beta_i\) so that it intersects
\(\beta_i\) transversely at 2 points. As described in \cite{OS:4d},
the ``Heegaard triple'' \((\Sigma, \pmb{\alpha}, \pmb{\beta},
\pmb{\gamma})\) defines a 4-dimensional cobordism \(X_{\pmb{\alpha}, \pmb{\beta},
\pmb{\gamma}}\) from \(M\sqcup \#^{G-1} (S^1\times S^2)\) to
\(M_K\). Associated to this  are holomorphic triangle maps \(F^\circ _{\pmb{\alpha}, \pmb{\beta},
\pmb{\gamma}}: HF^\circ( M)\otimes HF^-(\#^{G-1}(S^1\times S^2))\to HF^\circ (M_K)\) (Cf. e.g. \cite{OS:4d} Section 2.3).
The elementary cobordism \(W_K\) may be recovered from \(X_{\pmb{\alpha}, \pmb{\beta},
\pmb{\gamma}}\) by filling in the boundary component \(\#^{G-1}
(S^1\times S^2)\) with the boundary connected sum \(\#^{G-1}(S^1\times
B^3)\). The cobordism map \(F_{W_K}^\circ: HF^\circ( M)\to HF^\circ (M_K)\) is defined in terms of the
holomorphicm triangle maps \(F^\circ _{\pmb{\alpha},
  \pmb{\beta},\pmb{\gamma}}\) by evaluating the latter on the highest
degree generator (unique modulo signs) of \(HF^-(\#^{G-1}(S^1\times S^2))\simeq\bigotimes^{G-1}\hat{V}\). 

 Let \(\Delta\) denote a triangle. Extending
 the aforementioned re-interpretation of the Heegaard Floer chain
 maps, in \cite{Li} the Heegaard
 triple maps have an alternative definition as counting invariants of
 holomorphic curves in the 4-manifold \(\Delta\times \Sigma\), with
 boundary conditions specified by the triple \((\pmb{\alpha}, \pmb{\beta},
\pmb{\gamma}\)). This interpretation is not directly compatible with
our proof in
\cite{KLT} of \(HF=HM\) through intermediate Floer homologies . Instead, we reinterprete \(F_{W_K}\) as a counting
invariant of holomorphic sections in a 4-dimensional Lefschetz
fibration over \(\mathbb{R}\times I\) that contains one single
singular points. The boundary conditions are specified by
\((\pmb{\alpha}, \pmb{\beta})\), and in this picture, \(\gamma_G\)
arises as the Dehn twist of \(\beta_G\) along the vanishing cycle of
the Lefschetz singularity. The relation between this interpretation
and Lipshitz's has antecedent analogs, see e.g. \cite{S}. 

With this 4-dimensional interpretation of \(F_{W_K}^\circ\), we may now
modify the cobordism \(W_K\) in a straightforward manner to get an
auxiliary cobordism \(W_Y\) between \(Y_M\) and \(Y_{M_K}\). This
cobordism induces maps between the intermediate Floer homologies
\(H_{SW}^\circ\) of \(M\) and \(M_K\), and Theorems \ref{thm:A},
\ref{thm:B} have direct analogs that compute this map in terms of
the Heegaard Floer map \(F_{W_K}^\circ\). Meanwhile, an application of the
filtered connected sum formula, Theorem \ref{thm:conn_sum}, allows one to
relate this map with the Seiberg-Witten cobordism maps \(\mathring{HM}(W_K)\).

\paragraph{Acknowledgements.} This research is supported in part by
Hong Kong Research Grant Counsel GRF \# 401913.

\end{document}